\documentclass[11pt,reqno]{amsart}
\usepackage{amssymb,accents}

\def\noi{\noindent}

\newtheorem{Thm}{Theorem}[section]
\newtheorem{Def}[Thm]{Definition}
\newtheorem{Lm}[Thm]{Lemma}
\newtheorem{Prop}[Thm]{Proposition}
\newtheorem{Cor}[Thm]{Corollary}
\newtheorem{Rem}[Thm]{Remark}
\newtheorem{state}{Theorem}

\def\cal{\mathcal}
\def\Bbb{\mathbb}
\def\mf{\mathfrak}

\def\<{\langle}
\def\>{\rangle}

\def\a{\alpha}
\def\b{\beta}
\def\d{\delta}
\def\D{\Delta}
\def\th{\theta}
\def\eps{\varepsilon}
\def\l{\lambda}
\def\L{\Lambda}

\def\Re{\Bbb R}
\def\F{\Bbb F}
\def\C{\Bbb C}
\def\Z{\Bbb Z}

\def\R{\cal R}
\def\S{\cal S}
\def\P{\cal P}
\def\H{\cal H}
\def\O{\cal O}
\def\A{\cal A}

\def\G{\mf h}

\def\W{\ring W}
\def\w{\ring w}
\def\RR{\ring R}
\def\Q{\ring Q}
\def\h{\ring {\cal H}}
\def\Gc{\ring {\mf h}}
\def\E{\ring E}
\def\bo{\ring {\mf b}}

\begin{document}
\title[nonsymmetric macdonald polynomials]
{Nonsymmetric Macdonald polynomials and Demazure characters}
\author{Bogdan Ion} 
\address{
Department of Mathematics, 
Princeton University,
Princeton NJ-08544
}
\email{bogdan@math.princeton.edu}
\begin{abstract}
We establish a connection between a specialization of the nonsymmetric 
Macdonald
polynomials and the Demazure characters of the corresponding affine Kac-Moody
algebra. This allows us to obtain a representation-theoretical interpretation
of the coefficients of the expansion of the specialized
symmetric Macdonald polynomials in the basis
formed by the irreducible characters of the associated finite Lie algebra.
\end{abstract}
\maketitle

\thispagestyle{empty}
\section*{Introduction}
\bigskip
Generalizing the characters of compact simple Lie groups I.G. Macdonald
associated to each irreducible root system a family of orthogonal 
polynomials $\P_\l(q,t)$ indexed by anti-dominant weights  and 
which are invariant under the
action of the Weyl group. These polynomials depend rationally on 
parameters $q$ and $t=(t_s,t_l)$ and for particular values of these
parameters reduce to familiar objects in representation theory:
\begin{enumerate}
\item when $q=t_s=t_l$ they are
equal to $\chi_\l$ the Weyl characters
of the corresponding root system (in particular they are independent of $q$);
\item when $q=0$ they are the polynomials that give the values of zonal 
spherical functions on a semisimple $p$-adic Lie group relative to a maximal
compact subgroup;
\item when $t_s=q^{k_s}, \ t_l=q^{k_l}$ and $q$ tends to $1$ they are
the polynomials that give the values of zonal 
spherical functions on a real symmetric space $G/K$ that arise from finite 
dimensional spherical representations of $G$. Here $k_s$, $k_l$ are the 
multiplicities of the short, respectively long, restricted roots.
\end{enumerate}

The nonsymmetric Macdonald polynomials $E_\l(q,t)$ (indexed this time by the 
entire weight 
lattice) 
were first introduced 
by E. Opdam \cite{o} in the differential setting and then by I. Cherednik 
\cite{c2} in full
generality. 
Unlike the symmetric polynomials, their 
representation-theoretical meaning is still unexplored. 
At present time their main 
importance consists in the fact that they form the common spectrum of a 
family of commuting operators (the Cherednik operators) which play a 
preponderant role in the representation theory of affine Hecke algebras and 
related harmonic analysis.

\medskip
It became clear, especially from the work of Cherednik, that we can  
in fact construct  such families of polynomials for every irreducible
{\sl affine root system}. 
From this point of view, the objects studied by Cherednik 
(\cite{c1},\cite{c2}, \cite{c3})
are the 
polynomials attached to  {\sl reduced twisted affine root systems}, and the 
Koornwinder polynomials, studied by S. Sahi (\cite{s2},\cite{s3}), are the 
polynomials attached to {\sl non-reduced affine root systems}.

\medskip
 This paper was inspired by the result of 
Y. Sanderson \cite{sa} who established a connection between a 
specialized version of the nonsymmetric Macdonald polynomials 
($E_\l(q,\infty)$ in our notation)
and 
the characters of a certain Demazure modules  $E_{w_\l}(\tilde\l)$ 
of the irreducible affine Lie algebra
(see Section \ref{sect1} for the definitions of the ingredients)
in the case of an irreducible root system of type $A_n$. 
Extrapolating from \cite{sa} we 
establish here the same connection for all irreducible affine root systems
for which {\sl the affine simple root is short}. This condition identifies
precisely the polynomials studied by Cherednik and Sahi. The proof rely
heavily on the method of intertwiners in double affine Hecke algebras.

\begin{state}\label{teorema1}
For an affine root system as above and 
any weight $\l$ we have,
$$
E_\l(q,\infty)=q^{\left(\L_0,w_\l\<\tilde\l\>\right)}
\chi(E_{w_\l}(\tilde\l))\ .
$$
\end{state}
The remaining cases: $B_n^{(1)}$, $C_n^{(1)}$, $F_4^{(1)}$ and $G_2^{(1)}$
exhibit some special features. For example, the formula of the affine 
intertwiner as an element of the double affine Hecke algebra takes a different
form (see \cite{ion}). Computations suggest that the action 
on the weight lattice
of the 
degeneration of this  affine intertwiner at
$t=\infty$
does not equal the action of the affine Demazure operator, but a different 
action with similar properties.

The connection between nonsymmetric Macdonald polynomials and Demazure 
characters
allows a  representation-theoretical interpretation of the 
coefficients of the expansion 
of the symmetric polynomials in the basis formed by the irreducible
characters of the associated finite Lie algebra. Our second  result
is the following
\begin{state}\label{teorema2}
For an affine root system as above and 
any anti-dominant weight $\l$ the symmetric polynomial 
$P_\l(q,\infty)$ can be written as a sum
$$
P_\l(q,\infty)=
\sum_{\mu\leq\l}d_{\l\mu}(q)\chi_\mu
$$
where $d_{\l\mu}(q)$ is a
polynomial in $q^{-1}$ with positive integer coefficients.
\end{state}
Let us mention that in the $A_n$ case, as explained in \cite{sa}, the 
positivity of the above coefficients is closely related to the 
positivity of the Kostka-Foulkes polynomials via the duality of 
the two variable Kostka functions. Another consequence of the Theorem
\ref{teorema1} is the following
\begin{state}\label{teorema3}
For an affine root system as above and 
any weight $\l$ we have,
$$
E_\l(\infty,\infty)=
\chi(\E_{\w_\l w_\circ}(\l_+))\ .
$$
\end{state}
This relates the specialization of the nonsymmetric Macdonald polynomials
$$
E_\l(\infty,\infty)=\lim_{q\to \infty}\lim_{t\to \infty}E_\l(q,t)
$$
to the Demazure characters of the finite irreducible Lie algebras. The order in
which we compute the above limits seems to be irrelevant.

{\sl Acknowledgement.} I want to acknowledge my deep gratitude to Professor 
Siddhartha Sahi for his
generous and inspiring guidance.

\bigskip
\section{Preliminaries}\label{sect1}
\bigskip

\subsection{The affine Weyl group}\label{sect1.1}

For the most part we shall adhere to the notation in \cite{kac}.
Let $A=(a_{ij})_{0\leq i,j\leq n}$ be an irreducible \emph{affine} 
Cartan matrix, 
$S(A)$ the Dynkin diagram and  $(a_0,\dots, a_n)$  the numerical labels of 
$S(A)$ in Table Aff from \cite{kac}, p.48-49. We denote by 
$(a_0^\vee,\dots, a_n^\vee)$ the labels of the Dynkin diagram $S(A^t)$ of the 
dual algebra which is obtained from $S(A)$ by reversing the direction
of all arrows and keeping the same enumeration of the vertices. 
Let $({\G}, R, R^{\vee})$ be a realization of $A$
and let
$(\Gc, \RR, 
\RR^{\vee})$ 
be the associated finite root system (which is a realization of the Cartan
matrix $\ring A = (a_{ij})_{1\leq i,j\leq n}$). From 
this data
one can construct 
an {\sl affine Kac-Moody algebra} ${\mf g}$, 
respectively a {\sl finite
Lie algebra} $\ring {\mf g}$ such that $\G$, $\Gc$ become the 
corresponding Cartan subalgebras and $R$, $\RR$ become the corresponding root
systems. Note also that $\ring {\mf g}$ is a subalgebra of
${\mf g}$. We refer to \cite{kac} for the details of this construction.
If we denote by $\{\a_i\}_{0\leq i\leq n}$ a basis of $R$ 
such that $\{\a_i\}_{1\leq i\leq n}$ is a basis of $\RR$
we have the following description
$$
{\G}^*={\Gc^*} + {\mathbb R}\delta + 
{\mathbb R}{\Lambda}_0\ ,
$$
where $\d=\sum_{i=0}^n a_i\a_i$. The vector space ${\G}^*$ has a canonical
scalar product defined as follows
$$
(\a_i,\a_j):=d_i^{-1}a_{ij}\ ,\ \ \ \
(\L_0,\a_i):=\d_{i,0}a_0^{-1}\ \ \ \text{and}\ \ \ (\L_0,\L_0):=0,
$$
with $d_i:= a_ia_i^{{\vee}-1}$ and $\d_{i,0}$  Kronecker's delta.
As usual, $\{\a_i^\vee:=d_i\a_i\}_{0\leq i\leq n}$, $\{\l_i\}_{1\leq i\leq n}$
and $\{\l_i^\vee\}_{1\leq i\leq n}$ are the coroots, fundamental weights and
fundamental coweights. Denote by  
$P=\oplus_{i=1}^n\Z\l_i$ and $\Q=\oplus_{i=1}^n\Z\a_i$
the weight lattice, respectively the root lattice and let
$$
\rho:= \frac{1}{2}\sum_{\a\in \RR_+}\a^\vee=\sum_{i=1}^n\l_i^\vee.
$$

Given $\a\in R$, $x\in \G^*$ let
$$
s_\a(x):=x-\frac{2(x,\a)}{(\a,\a)}\a\ .
$$
The {\sl affine Weyl group} $W$ is generated by all $s_\a$ (the simple
reflexions $s_i=s_{\a_i}$ are enough). 
The {\sl finite Weyl group} $\W$ is
the subgroup generated by $s_1,\dots,s_n$. An important role is played by
$\th=\delta -a_0\a_0$. Remark that $a_0=1$ in all cases except for 
$A=A_{2n}^{(2)}$, when $a_0=2$.
For $s$ a real number, $\G^*_s=\{ x\in\G\ ;\ (x,\d)=s\}$ is the 
level $s$ of $\G^*$. We have 
$$
\G^*_s=\G^*_0+s\L_0=\Gc^*+{\Bbb R}\d+s\L_0\ .
$$
 The action of $W$ 
preserves each of the $\G^*_s$ and  we can identify each of the $\G^*_s$
canonically with $\G^*_0$ and obtain an (affine) action of $W$ on $\G^*_0$.
If $s_i\in W$ is a simple reflexion, write $s_i(\cdot)$ for the regular
action of $s_i$ on $\G^*_0$
and $s_i\<\cdot\>$ for the affine action of 
$s_i$ on $\G^*_0$ corresponding to the level one action. These actions differ
only for $s_0$\ :
\begin{eqnarray*}
s_0(x) & = & s_\th(x)+(x,\th)\d\ ,\\
s_0\<x\> & = & s_\th(x)+a_0^{-1}\th+(x,\th)\d-a_0^{-1}\d\ .
\end{eqnarray*}
By $s_i\cdot$ we denote the affine action of $W$ on $\Gc^*$
$$
s_0\cdot x = s_\th(x)+a_0^{-1}\th\ .
$$

We will be interested in the cases when $\a_0$ is a {\sl short} root. This 
happens precisely when the affine root system is {\sl twisted} or
{\sl simply laced untwisted}. Under these conditions
define the fundamental alcove as 
$$
\A:=\{ x\in \Gc^*\ | \ 
(x+\L_0,\a_i^\vee)\geq 0\ ,\ 0\leq i\leq n\}.
$$ 
The non-zero elements of
$\O=P\bigcap \A$ are the so-called minuscule weights. Let us 
remark that the orbits of the affine action of $W$ on 
$P$ contains
a unique $\l_i\in \A$ (to keep the notation consistent we set $\l_0=0$).
\bigskip

\begin{quotation}{\sl
In all what follows we assume our affine root system to be such that
the affine simple root $\a_0$ is short (this condition includes of course the
case when all roots have the same length).}
\end{quotation}

\medskip
\subsection{The Bruhat order}
Let us first establish some notation.
For each $w$ in $W$ let $l(w)$ be the length of a 
reduced (i.e. shortest) decomposition of $w$ in terms of the $s_i$.
We have
$
l(w)=|\Pi(w)|
$
where
$
\Pi(w)=\{\a\in R_+\ |\ w(\a)\in R_-\}\ .
$
If $w=s_{j_p}\cdots s_{j_1}$ is a reduced decomposition, then
$$
\Pi(w)=\{\a^{(i)}\ |\ 1\leq i\leq p\},
$$
with $\a^{(i)}=s_{j_1}\cdots s_{j_{i-1}}(\a_{j_i})$.
For each weight $\l$ define 
$\l_-$, respectively $\tilde\l$,
to be the unique element in 
$\W\hspace{-0.1cm}\l$, respectively $W\cdot\l$,
which is an anti-dominant weight, respectively an element of $\O$
(that is a minuscule weight or zero), and 
$\w_\l\hspace{-0.2cm}^{-1}\in\ \W$, $w_\l^{-1}\in W$,  
to be the unique minimal length elements 
by which this is achieved. Also, for each weight $\l$ define $\l_+$ to be the 
unique element in $\W\l$ which is dominant and denote by $w_\circ$ 
the maximal 
length element in $\W$.
\begin{Lm}\label{lema1} 
With the notation above, we have
\begin{enumerate}
\item[i)]
$
\Pi(\w_\l\hspace{-0.2cm}^{-1})=\{\a\in \RR_+\ |\ (\l,\a)>0 \}\ ;
$
\item[ii)]
$
\Pi(w_\l^{-1})=\{\a\in R_+\ |\ (\l+\L_0,\a)<0 \}\ .
$
\end{enumerate}
\end{Lm}
\begin{proof}
Straightforward. See
Theorem 1.4 of \cite{c1} for a full argument.
\end{proof}

The Bruhat order is a partial order on any Coxeter group.
For its basic properties see
Chapter 5 in \cite{humph}. Let us list a few of them 
(the first two properties completely characterize the Bruhat order):
\begin{enumerate}
\item For each $\a\in R_+$ we have $s_\a w<w$ iff $\a$ is in $\Pi(w^{-1})$ ;
\item $w'< w$ iff $w'$ can be obtained by omitting some factors
in a fixed reduced decomposition of $w$ ;
\item if $w' \leq w$ then either $s_i w' \leq w$ or 
$s_i w' \leq s_iw$ (or both). 
\end{enumerate}
We can use the Bruhat order on $W$ do define a partial order on the weight 
lattice: 
if $\l,\mu\in P$ then by definition $\l<\mu$ iff $w_\l<w_\mu$.
\begin{Lm}\label{lemma2}
Let $\l$ be a weight such that $s_i\cdot \l\neq \l$ for some  
$0\leq i\leq n$. Then
$w_{s_i\cdot \l}=s_iw_\l$.
\end{Lm}
\begin{proof}
Because $l(s_iw_\l)=l(w_\l)\pm 1$ and 
$l(s_iw_{s_i\cdot \l})=l(w_{s_i\cdot\l})\pm 1$ we have  four possible 
situations depending on the choice of the signs in the above relations. The
choice of a plus sign in both relations translates in 
$\a_i\not\in\Pi(w_\l^{-1})$ and 
$\a_i\not\in\Pi(w_{s_i\cdot\l}^{-1})$ which by Lemma \ref{lema1} and our 
hypothesis implies that $(\a_i,\l+\L_0)>0$ and 
$(\a_i,s_i\cdot\l+\L_0)> 0$ (contradiction). The same argument shows that
the choice of a minus sign in both relations is impossible. Now, we can 
suppose 
that $l(s_iw_\l)=l(w_\l)+ 1$ and 
$l(s_iw_{s_i\cdot \l})=l(w_{s_i\cdot\l})- 1$, the other case being treated 
similarly. Using the minimal length properties of $w_\l$ and $w_{s_i\cdot\l}$
we can write
$$
l(w_\l)+ 1=l(s_iw_\l)\geq l(w_{s_i\cdot\l})=
l(s_iw_{s_i\cdot \l})+1\geq l(w_\l)+ 1
$$
which shows that $l(s_iw_\l)=l(w_{s_i\cdot\l})$. Our conclusion follows from
the uniqueness of $w_{s_i\cdot\l}$.
\end{proof}
An immediate consequence is the following
\begin{Lm}\label{lemma1.3}
Let $\l$ be a weight such that $s_i\cdot \l\neq \l$ for some  
$0\leq i\leq n$. Then
$s_i\cdot\l>\l$ iff $(\a_i,\l+\L_0)> 0$.
\end{Lm}
\begin{Lm}\label{lemma3}
Let $\l$ be a weight such that $s_i\cdot \l\neq \l$ for some  
$0\leq i\leq n$. Then
$\w_{s_i\cdot \l}=s_i\w_\l$ if $i\neq 0$ and $\w_{s_0\cdot \l}=s_\th\w_\l$.
\end{Lm}
\begin{proof}
We can prove the statement for $i\neq 0$ with the same arguments as in 
Lemma \ref{lemma2}. The remaining statement was essentially proved in 
Lemma 3.3 of \cite{s4}.
\end{proof}
\begin{Def}
If $\l$ and $\mu$ are weights such that $\l-\mu\in\Q$, we say that 
the weight $\nu$ is a convex combination of $\l$ and $\mu$ if 
$\nu=(1-\tau)\l+\tau\mu$ such that $0\leq\tau \leq 1$ and $\l-\nu\in\Q$.
\end{Def}
The following result was proved in Lemma 5.5 of \cite{s3} for a particular
affine Weyl group, but the proof provided there works in general.
\begin{Lm}\label{lemma1.6}
Let $\l$ be a weight such that $s_i\cdot \l\geq \l$ for some  
$0\leq i\leq n$. If $\nu$ is a convex combination of $\l$ and $s_i\cdot\l$,
then
$\nu<\l<s_i\cdot\l$.
\end{Lm}
For each weight $\l$ define $\overline\l=\l+\w_\l(\rho)$. As a consequence of
Lemma \ref{lemma3} we have the following 
\begin{Prop}\label{prop1.7}
Let $\l$ be a weight such that $s_i\cdot \l\neq \l$. Then 
$$
s_i\cdot\overline\l=\overline{s_i\cdot\l}\ \ .
$$
\end{Prop}
\medskip
\subsection{Demazure modules characters}
Recall that ${\mf g}$ is the Kac-Moody affine Lie algebra associated with the 
irreducible affine Cartan matrix $A$. For the results in this section we 
refer to \cite{kum}.
Let $\L$ be an integral {\sl dominant}
weight of ${\mf g}$ 
and let $V=V(\L)$ be the unique irreducible highest weight  
${\mf g}$-module with highest weight $\L$. For each $w\in W$ the weight space
$V_{w(\L)}$ is one-dimensional. Consider $E_w(\L)$, the 
${\mf b}$-module generated by $V_{w(\L)}$, where ${\mf b}$ is the the Borel 
subalgebra of ${\mf g}$. 
The  $E_w(\L)$, called the Demazure modules, are finite
dimensional vector spaces. If $\l$ is an element of $\O$, then $\l+\L_0$
is a dominant weight. In such a case we write $E_w(\l)$ for $E_w(\l+\L_0)$.
To a Demazure module $E_w(\L)$ we can associate its
character
$$
\chi(E_w(\L))=\sum_{\Upsilon\text{ weight }} \text{dim}(E_w(\L)_{\Upsilon})
\cdot e^{\Upsilon} 
$$
which can be regarded as an element in 
$\P:=\C[q,q^{-1}][e^{\mu};\ \mu\in P]$ after we ignore the $e^{\L_0}$ factor 
and after we set $q=e^{-\d}$.
\begin{Def}
For each $0\leq i\leq n$ define an operator $\D_i$ acting on $\P$
$$
\D_ie^{\L}=
\frac{ e^{\L} -e^{-\a_i}e^{s_i\<\L\>}}{1-e^{-\a_i}}\ .
$$
\end{Def}
Let $w=s_{i_1}\cdots s_{i_j}$ be  a reduced decomposition. Then, we can define
$\D_w:=\D_{i_1}\cdots \D_{i_j}$ (the definition of $\D_w$ does not depend
on the choice of the reduced decomposition).  
\begin{Thm}\label{T1.9} Let $\l$ be an element of $\O$. Then
$$\chi(E_w(\l))=\D_w(e^\l)\ .$$ 
\end{Thm}
The above Theorem is a special case of the Demazure character formula for 
Kac-Moody algebras, which has proved in full generality by S. Kumar and
independently by O. Mathieu. We refer to Theorem 3.4 of \cite{kum} for the 
proof. 
The construction of the Demazure modules $\E_{\w}(\l)$
for the Lie algebra $\ring {\mf g}$
is completely  analogous (the role of ${\mf b}$ is played here by $\bo$ the 
Borel subalgebra of $\ring {\mf g}$).

\bigskip
\section{Nonsymmetric Macdonald polynomials}
\bigskip
In what follows we consider our root system to be {\sl reduced}. 
Recall that in this case $a_0=1$ and $\th$ is the highest {\sl short} root
of the associated finite root system.
The case of a non-reduced root system will be treated
in Section \ref{nonred}.

\subsection{The double affine Hecke algebra}
We introduce a field $\F$ (of parameters) as follows:
fix indeterminates $q$ and $t_0,\dots,t_n$ such that
$
t_i=t_j \text{ iff } d_i=d_j;
$
let $m$ be the lowest common denominator of the rational numbers
$\{(\a_i,\l_j)\ |\ 1\leq i,j\leq n \}$, and let $\F$ denote the field of 
rational functions in $q^{1/m}$ and $t_i^{1/2}$. Because in our case there
are at most two different root lengths we will also use the notation 
$t_l$, $t_s$ for $t_i$ if the corresponding simple root is long,
respectively short. 
 The algebra $\R=\F[e^\l;\l\in P]$ is the group 
$\F$-algebra of the lattice
$P$ and $\S$ is the subalgebra of $\R$ consisting of 
elements invariant under the finite Weyl group. 
For further use we also introduce the following 
group $\F$-algebras of the root lattice:
$\R_Y:=\F[Y_\mu;\mu\in \Q]$ and $\R_X:=\F[X_\b;\b\in \Q]$.
$\S_Y$ is the subalgebra of $\R_Y$ consisting of elements invariant under
the finite Weyl group.
\begin{Def}
The affine Hecke algebra $\H$ is the $\F$-algebra generated by elements
$T_0,\dots,T_n$ with relations
\begin{enumerate}
\item[(i)] the $T_i$ satisfy the same braid relations as the $s_i$;
\item[(ii)] for $0 \leq i\leq n$ we have
$$
T_i^2=(t_i^{1/2} -t_i^{-1/2})T_i +1.
$$
\end{enumerate}
\end{Def}
The elements $T_1, \dots, T_n$ generate the  {\sl finite Hecke algebra}
$\h$. There are natural bases of $\H$ and $\h$: 
$\{T_w\}_w$ indexed by
$w$ in $W$ and in
$\W$ respectively, where $T_w=T_{i_l}\cdots T_{i_1}$ if 
$w=s_{i_l}\cdots s_{i_1}$ is a reduced  expression of $w$
in terms of simple reflexions. 
There is another important description of the affine Hecke algebra due to
Lusztig \cite{l2}.
\begin{Prop}\label{prop2.2}
The affine Hecke algebra $\H$ is generated by the finite Hecke algebra 
and the group algebra $\R_Y$ such that the following relation is satisfied
for any
$\mu$ in the root lattice and
any $1\leq i\leq n$\ :
$$
Y_\mu T_i-T_iY_{s_i(\mu)} = (t_i^{1/2}-t_i^{-1/2})
\frac{Y_\mu-Y_{s_i(\mu)}}{1-Y_{\a_i}}\ .
$$
\end{Prop}
\begin{Rem}
In this description $T_0^{-1}=Y_\th T_{s_\th}$.
\end{Rem} 
Following Macdonald \cite{mac}, we call 
the family of commuting operators 
$\R_Y\subset \H$, {\sl Cherednik operators}.
\noi In order to state the next result we need the following notations:
for 
$\mu,\b\in \Q$ and $k\in \Z$,
$X_{\b+k\d}:=q^{-k}X_\b$ and  $Y_{\mu+k\d}:=q^kY_{\mu}$. 
For the next results we refer to 
Cherednik \cite{c1}\cite{c3}. 
\begin{Def}
The double affine Hecke algebra $\H^d$ is the $\F$-algebra generated by 
the affine Hecke algebra $\H$ and the group algebra $\R_X$ 
such that the following relation is satisfied for any
$\b$ in the root lattice and
any $0\leq i\leq n$:
$$
T_iX_\b-X_{s_i(\b)}T_i = (t_i^{1/2}-t_i^{-1/2})
\frac{X_\b-X_{s_i(\b)}}{1-X_{-\a_i}}\ .
$$
\end{Def}

\noi The following formulas define a faithful 
representation of $\H^d$ on $\R$ 
$$
\pi(T_i)e^\l=t_i^{1/2}e^{s_i(\l)} +(t_i^{1/2}-t_i^{-1/2})
\frac{e^\l-e^{s_i(\l)}}{1-e^{-\a_i}}\ ,\ \ \ 0\leq i\leq n
$$
$$
\pi(X_\b)e^\l=e^{\l+\b}\ ,\ \ \ \b\in \ \Q.
$$
\begin{Thm}\label{thm2.5}
Define $T_{\<0\>}=T_0^{-1}X_{\a_0}$. Then for all $\mu\in \Q$ and all $\l\in P$
$$
Y_\mu T_{\<0\>}-T_{\<0\>}Y_{s_0(\mu)} = (t_0^{1/2}-t_0^{-1/2})
\frac{Y_\mu-Y_{s_0(\mu)}}{1-Y_{\a_0}}\ ,
$$
$$
\pi(T_{\<0\>})e^\l=t_0^{1/2}e^{s_0\<\l\>} +(t_0^{1/2}-t_0^{-1/2})
\frac{e^\l-e^{s_0\<\l\>}}{1-e^{-\a_0}}\ .
$$
\end{Thm}

The irreducible affine root systems for which 
the affine simple root is long are $B_n^{(1)}$, $C_n^{(1)}$, $F_4^{(1)}$ 
and $G_2^{(1)}$. For these root systems
the formula of the element of the  
double affine Hecke algebra which plays the same role as $T_{\<0\>}$
takes a different
form (see \cite{ion}) which makes the computation of its action on $\R$
more difficult.

To avoid cumbersome notation we set  $T_{\<i\>}=T_i$ for 
$i\neq 0$ and write
$He^\l$ in place of $\pi(H)e^\l$ for any $H\in \H^d$.
 The following theorem follows directly from Lemma \ref{lemma1.6}.
\begin{Thm}\label{triang}
Suppose that
$\l\leq\gamma\leq s_i\cdot\gamma$ for some weights $\l$, $\gamma$ and 
$0\leq i\leq n$. 
Then
$$
T_{\<i\>}e^\l=t_i^{1/2}e^{s_i\<\l\>}\ +\ \text{lower terms}
$$
where by lower terms we mean a combination of $e^\b$ with $\b<s_i\cdot \l$.
\end{Thm}
\medskip
\subsection{Macdonald polynomials}\label{sect2.2}
Cherednik defined a certain scalar product on $\R$ (see \cite{c1} for details)
for which all operators in $\H$ became unitary operators.
In particular the adjoint of $Y_\mu$ is
$Y_{-\mu}$.
By 
${\bf q}^{(\mu+k\d,\overline\l)}$ we denote the element of $\F$
$$
q^{k+(\mu,\l)}\prod_{i=1}^{n}t_i^{-(\mu, \w_\l(\l_i^\vee))}\ .
$$

\noi For each $\l\in P$ we can construct
a $\F$-algebra morphism ${\rm ev}(\l):\R_Y\to \F$, which sends  
$Y_\mu$ to ${\bf q}^{(\mu,\overline \l)}$.
If $f$ is an element of $\R_Y$ we will write $f(\l)$ for ${\rm ev}(\l)(f)$.
Macdonald defined a basis $\{P_\l(q,t)\}$
of $\S$ 
which is indexed by {\sl anti-dominant} 
weights and which is completely characterized by the equations
\begin{equation}\label{mach}
f\cdot P_\l=f(\l)P_\l
\end{equation}
for any $f\in \S_Y$, 
and the condition that the coefficient of $e^{\l}$ in $P_\l(q,t)$ is $1$. 
The elements of this basis are called {\sl symmetric Macdonald polynomials}.

\medskip
Recently, a nonsymmetric version of the Macdonald polynomials was introduced
by Opdam \cite{o} in the differential case, Macdonald \cite{mac} (for
$t_i=q^k$, $k\in \Z_+$) and by Cherednik \cite{c2} in the general (reduced)
case and some of their properties were studied. 
For each weight $\l$ there is an unique element $E_\l(q,t)\in \R$ 
satisfying
the conditions
\begin{eqnarray}\label{scal}
        E_\l &=& e^\l\ + \text{ lower terms};\\
(E_\l,e^\mu) &=& 0 \text{\ \ \  for all\ \  } \mu<\l\ .\label{ss}
\end{eqnarray}
They form a $\F$-basis of $\R$ and they are the common eigenfunctions
of the Cherednik operators. In what follows we will find an explicit recursion
formula for the nonsymmetric Macdonald polynomials. In the course of doing
that we will give a more transparent proof of their existence and uniqueness.

\medskip
For all $0\leq i\leq n$ let us introduce the following elements of $\H^d$
called {intertwiners}
$$
I_i:=T_{\<i\>}(1-Y_{\a_i})-(t_i^{1/2}-t_i^{1/2})\ .
$$
The intertwiners where first introduced by Knop and Sahi 
\cite{knop},\cite{ks},\cite{s1} for ${\rm GL}_n$ and then by Cherednik
\cite{c3} in the general (reduced) case.
Their importance is the following: for any $\mu$ in the root lattice we have
\begin{equation}\label{intert}
Y_{\mu}I_i=I_iY_{s_i(\mu)}.
\end{equation}
This easily follows from Proposition \ref{prop2.2} and Theorem \ref{thm2.5}.
The next results can be proved following closely the ideas in \cite{s2} where
the non-reduced case was considered. 
For every weight $\l$ define 
$$
\R_{\l}=\{f\in\R\ | \ Y_\mu f={\bf q}^{(\mu,\overline \l)}f \ \text{for any } 
\mu\in M \}.
$$
\begin{Thm}\label{T32}
Let $\l$ be a weight such that $s_i\cdot \l\neq \l$. Then 
$I_i:\R_\l \to \R_{s_i\cdot \l}$ is a linear isomorphism.
\end{Thm}
\begin{proof}
Let $f$ be any element of $\R_\l$. 
Using the intertwining relation (\ref{intert}) and the Proposition 
\ref{prop1.7} 
we get
$$
Y_\mu(I_if)={\bf q}^{(\mu,\overline{s_i\cdot \l})}I_if\ .
$$
Therefore, $I_if$ is an element of $\R_{s_i\cdot \l}$. 
A short computation shows that
$$
I_i^2=t_i+t_i^{-1}-(Y_{\a_i}+Y_{-\a_i}),
$$
therefore $I_i^2$ acts as a constant on $\R_\l$. It is easy to see that our
hypothesis implies that this constant is nonzero, showing that $I_i^2$ and
consequently $I_i$ is an isomorphism.
\end{proof}
\begin{Thm}
The spaces $\R_\l$ are one-dimensional.
\end{Thm}
\begin{proof}
The proof is very similar with the proof of the corresponding result
(Theorem 6.1) in \cite{s2}. The only difference is that we have to use 
the fact that $\O$ is a set of
representatives for the orbits of the affine action of $W$ on $P$, and 
the fact that  $e^\l$ is in
$\R_\l$ for  $\l\in \O$. From the proof also follows that an element
in $\R_\l$ is uniquely determined by the coefficient of $e^\l$ in $f$.
\end{proof}
This result makes possible the following definition.
\begin{Def}
For any weight $\l$ define
the nonsymmetric Macdonald polynomial $E_\l(q,t)$ to be the unique 
element in $\R_\l$ in which the coefficient of $e^\l$ is $1$. If $k\in \Z$
then denote $E_{\l+k\d}(q,t)=q^{-k}E_\l(q,t)$.
\end{Def}
For each anti-dominant weight $\l$ we write $\R^\l$ for the subspace of $\R$
spanned by  $\{E_\mu\ |\ \mu\in\ \W\hspace{-0.1cm}\l \}$.
The connection with the symmetric Macdonald polynomials is the  following.
\begin{Cor}
The polynomial $P_\l(q,t)$ can be characterized as the 
unique 
$\W$- invariant element in $\R^\l$ for which the coefficient of $e^\l$ equals
$1$.
\end{Cor}
\begin{proof}
The result follows from the characterization (\ref{mach}).
\end{proof}
\begin{Def}
Let $C$ be the element of the finite Hecke algebra defined by 
$C:=(\sum_{w\in\W}\chi(T_w)^2)^{-1}\sum_{w\in\W}\chi(T_w)T_w$, where $\chi$
is the one dimensional representation of $\h$ defined by $\chi(T_i)=t_i^{1/2}$.
\end{Def}
\begin{Cor}\label{38}
$\pi(C)$ is a projection from $\R^\l$ to $\F P_\l$.
\end{Cor}
\begin{proof}
An easy calculation as in Lemma 2.5 of \cite{s1} shows that 
$T_iC=t_i^{1/2}C$ for any $1\leq i\leq n$, hence $T_i(Cf)=t_i^{1/2}Cf$ for
all $f\in \R$. This implies that $Cf$ is $\W$-invariant, and so it must be a
multiple of $P_\l$. Moreover, $C$ acts as identity on $\S$.
\end{proof}
For any weight $\l$ and any $0\leq i\leq n$ define the operator 
$G_{i,\l}(q,t)$
as follows
$$
G_{i,\l}:=t_i^{-1/2}T_{\<i\>}\ \ \ \text{if} \ (\l+\L_0,\a_i)=0\ ,
\ \ \ \text{and}
$$
$$
G_{i,\l}:= (1-{\bf q}^{-(\a_i,\overline\l)})t_i^{-1/2}T_{\<i\>}+
{\bf q}^{-(\a_i,\overline\l)}(1-t_i^{-1})\ \ \ 
\text{if} \ (\l+\L_0,\a_i)\neq0\ .
$$
\begin{Thm}\label{T39}
Let $\l$ be a weight such that $(\l+\L_0,\a_i)\geq 0$. Then
\begin{equation}\label{recursion}
G_{i,\l}E_\l=
(1-{\bf q}^{-(\a_i,\overline\l)})E_{s_i\<\l\>}\ .
\end{equation}
\end{Thm}
\begin{proof}
When $(\l+\L_0,\a_i)=0$ the statement follows straightforward from  
(\ref{scal}), 
(\ref{ss}) and
from the Theorem \ref{triang}. For the remaining case,
using Theorem \ref{T32}, all we need is to compute the coefficient
of $e^{s_i\<\l\>}$ in $G_{i,\l}E_\l$ which by Theorem \ref{triang} can be
 shown 
to be $(1-{\bf q}^{-(\a_i,\overline\l)})$.
\end{proof}
\medskip
\subsection{The specialization at $t=\infty$}\label{sect2.3}
Our goal is to define the specialization of the polynomials $E_\l(q,t)$ at 
$t=\infty$ (that means $t^{-1}=0$)
and to obtain recursion formulas for them as in Theorem \ref{T39}. 
In order to do this we have to closely examine the coefficients of the
$E_\l$ and make sure that their limit exists. In fact, we can
suitably re-normalize the $E_\l$ such that all the coefficients in this
re-normalization are polynomials in $t_i^{-1}$ and the normalizing
factor approaches $1$ when $t$ tends to infinity. This will show that 
the limit of each of the coefficients of the $E_\l$ exists and it is bounded.

Recall 
$w_\l$ be the unique minimal length element of $W$ such that
$w_\l\cdot \tilde\l=\l$. Let 
$w_\l=s_{j_l}\cdots s_{j_1}$ be a reduced decomposition. Then,
\begin{equation}\label{1}
\Pi(w_\l)=\{\a^{(i)}:=
s_{j_1}\cdots s_{j_{i-1}}(\a_{j_i})\ |\  
1\leq i\leq l\}\ .
\end{equation}
This means in particular that $\a^{(j)}\in R_+$ and 
$w_\l(\a^{(j)})\in R_-$. 
Define 
\begin{equation}\label{1b}
\l_{(i)}:=s_{j_{i-1}}\cdots s_{j_1}\cdot \tilde\l\ ,
\end{equation}
for any $1\leq i\leq l+1$. Therefore, $\l_{(1)}=\tilde\l$ and 
$\l_{(l+1)}=\l$. The 
key property of the $\l_{(i)}$ is that
\begin{equation}\label{2}
(\l_{(i)}+\L_0,\a_{j_i})>0\ .
\end{equation}
This easily follows from (\ref{1}). Moreover, (\ref{2}) implies that
$\a_{j_i}\in \Pi(\w_{\l_{(i)}}\hspace{-0.5cm}^{-1}\, )$ if 
$j_i\neq 0$, meaning that 
$\w_{\l_{(i)}}\hspace{-0.5cm}^{-1}(\a_{j_i})$ is in $\RR_-$, respectively that
$\th\not\in \Pi(\w_{\l_{(i)}}\hspace{-0.5cm}^{-1}\, )$ 
if $j_i=0$,  meaning that
$\w_{\l_{(i)}}\hspace{-0.5cm}^{-1}(\th)$ is in $\RR_+$. 

\medskip
\noi Now, for all $1 \leq j\leq l$,\ \  all the exponents 
in the monomial
${\bf q}^{(\a_{j_i},\overline \l_{(i)})}$  
are positive integers and at least one of exponents the $t_i$ is nonzero. 
Define the re-normalization of $E_\l(q,t)$ to be 
$$
\prod_{i=1}^l(1-{\bf q}^{-(\a_{j_i},\overline \l_{(i)})})E_\l(q,t)\ .
$$
This formula (modulo a $q$ factor) is obtained by applying the recursion 
formula (\ref{recursion}) successively, starting with
$e^{\tilde\l}$. From this description it is clear that the 
powers of the $t_i$ appearing the expansion of  this re-normalization of 
$E_\l(q,t)$
are all {\sl negative} and therefore our desired specialization at $t=\infty$
is well defined. We denote by  $E_\l(q,\infty)$ this specialization.
This re-normalization
does not depend on the choice of the reduced decomposition of $w_\l$.
Remark also that the coefficient of $e^\l$ in $E_\l(q,\infty)$ is $1$. 
For each anti-dominant weight $\l$ we write $\R^\l(\infty)$ 
for the linear subspace 
spanned by  $\{E_\mu(q,\infty)\ |\ \mu\in\W\hspace{-0.1cm}\l \}$.
The polynomial $P_\l(q,\infty)$ is defined to be the  
unique 
$\W$- invariant element in $\R^\l(\infty)$ for which the coefficient of 
$e^\l$ equals
$1$.
\bigskip
\section{Nonsymmetric Koornwinder polynomials}\label{nonred}

In this section we will consider the case of a non-reduced root system.
Recall that in this case $A=A^2_{(2n)}$, $a_0=2$, $\th$ is the highest
root and $\O=\{0\}$.
\subsection{The recursion relation} 
The results in this section are due to 
Sahi \cite{s2}, \cite{s3}.
We introduce the field $\F$ as follows: fix indeterminates
$q$, $u=(u_0, u_n)$ and $t_0,\cdots,t_n$ identified as before; the field
$\F$ is  the field of rational functions in their square roots. 
We also define
$$a =t_n^{1/2}u_n^{1/2},\ \ 
b=-t_n^{1/2}u_n^{-1/2}, \ \ 
c=q^{1/2}t_0^{1/2}u_0^{1/2},\ \  
d= -q^{1/2}t_0^{1/2}u_0^{-1/2}.
$$
Note that in this case we have three different root lengths, therefore
$t=(t_s,t_m,t_l)$, where $t_s=t_0$, $t_l=t_n$ and $t_m=t_i$ for any 
$i\neq 0,n$.
As before $\R=\F[e^\l;\l\in P]$ is the group 
$\F$-algebra of the lattice
$P$ and $\S$ is the subalgebra of $\R$ consisting of 
elements invariant under the finite Weyl group. 
Also, define
$\R_Y:=\F[Y_\mu;\mu\in P]$ and $\R_X:=\F[X_\b;\b\in P]$.
$\S_Y$ is the subalgebra of $\R_Y$ consisting of elements invariant under
the finite Weyl group. The lattice $P$ can be identified with $\Z^n$ such
that the scalar product we defined in Section \ref{sect1.1}
is the canonical scalar product on
$\Re^n$. If
$\eps_1,\cdots,\eps_n$ are the unit vectors in $\Z^n$, then our choice of the
basis for the affine root system is
$$
\a_0=\frac{1}{2}\d+\eps_1, \  \ \a_i=-\eps_i+\eps_{i+1},\ \ \a_n=-2\eps_n\ .
$$

\medskip
\noi The double affine Hecke algebra in this case has a more complicate 
description (see \cite{s2} for details). We describe here 
only its action on $\R$:
\begin{itemize}
\item
$
T_0e^\l:= t_0^{1/2}e^\l+t_0^{-1/2} 
\frac{
(1-c e^{-\eps_1})(1-de^{-\eps_1})
}
{1-qe^{-2\eps_1} 
}(e^{s_0(\l)}-e^\l)
\ ,
$

\item
$
T_{\<0\>}e^{\l}:=T_0^{-1}e^{\l+\a_0}
\ ,
$
\item
$
T_{\<i\>}e^\l=T_ie^\l := t_i^{1/2}e^\l+t_i^{-1/2} 
\frac{(1- t_i e^{-\a_i})}{
(1-e^{-\a_i})} (e^{s_i(\l)}-e^\l)\ ,\quad i\neq 0,n,
$
\item
$
T_{\<n\>}e^\l=T_ne^\l:= t_n^{1/2}e^\l+t_n^{-1/2} 
\frac{
(1-a e^{\eps_n})(1-be^{\eps_n})
}
{1-qe^{2\eps_n} 
}(e^{s_n(\l)}-e^\l)
\ .
$
\end{itemize}
The commutative algebra $\R_Y$ embeds in the Hecke algebra
as follows
$$
Y_{\eps_i}= (T_i\cdots T_{n-1})(T_n\cdots T_0)
(T_1^{-1}\cdots T_{i-1}^{-1})\ .$$
The action of $\R_Y$ can be simultaneously diagonalized and the 
{\sl nonsymmetric
Koornwinder polynomials} $E_\l(q,t,u)$ are the corresponding eigenbasis. The
eigenvalues are given as follows:
by 
${\bf q}^{(\mu+k\d,\overline\l)}$ we denote the element of $\F$
$$
q^{k+(\mu,\l)}(t_0t_n)^{-(\mu, \w_\l(\l_n^\vee))}
\prod_{i=1}^{n-1}t_i^{-(\mu, \w_\l(\l_i^\vee))}\ .
$$

\noi For each $\l\in P$ we can construct
a $\F$-algebra morphism ${\rm ev}(\l):\R_Y\to \F$, which sends  
$Y_\mu$ to ${\bf q}^{(\mu,\overline \l)}$.
If $f$ is an element of $\R_Y$ we will write $f(\l)$ for ${\rm ev}(\l)(f)$.
The {\sl symmetric Koornwinder polynomials} $\{P_\l(q,t,u)\}$ form a basis
of $\S$ 
which is indexed by {\sl anti-dominant} 
weights. They are  completely characterized by the equations
\begin{equation}
f\cdot P_\l=f(\l)P_\l
\end{equation}
for any $f\in \S_Y$, 
and the condition that the coefficient of $e^{\l}$ in $P_\l(q,t,u)$ equals 
$1$. In the same manner as is Section \ref{sect2.2} we define 
for any weight $\l$
the vector spaces
$\R_\l$ and $\R^\l$. 
\begin{Prop}
The polynomial $P_\l(q,t,u)$ can be characterized as the 
unique 
$\W$- invariant element in $\R^\l$ for which the coefficient of $e^\l$ equals
$1$.
\end{Prop}
\noi For any weight $\l$ and any $0\leq i\leq n$ such that $(\l+\L_0,\a_i)=0$
define the operator $G_{i,\l}(q,t)$
as follows
$$
G_{i,\l}:=t_i^{-1/2}T_{\<i\>}\ .
$$
If $(\l+\L_0,\a_i)\neq 0$ we define 
$$
G_{i,\l}:= (1-{\bf q}^{-(\a_i,\overline\l)})t_i^{-1/2}T_{i}+
{\bf q}^{-(\a_i,\overline\l)}(1-t_i^{-1})\ \ \ 
\text{for } i\neq 0\ \ \ \text{and}
$$
$$
G_{0,\l}:= t_0^{-1/2}((1-{\bf q}^{-(\d-\th,\overline\l)})T_{\<0\>} 
+{\bf q}^{-(\a_0,\overline\l)}(u_n^{1/2}-u_n^{-1/2})
+(u_0^{1/2}-u_0^{-1/2}))\ .
$$
\begin{Thm}\label{T39koor}
Let $\l$ be a weight such that $(\l+\L_0,\a_i)\geq 0$. Then
\begin{equation}\label{recursionkoor}
G_{i,\l}E_\l=
(1-{\bf q}^{-(\a_i,\overline\l)-\d_{i,0}(\a_0,\overline\l)} )
E_{s_i\<\l\>}\ .
\end{equation}
\end{Thm}
\medskip
\subsection{The specialization at $u=(t_0,1)$,  $t=\infty$}

First, there is of course no problem in specializing $u_0:=t_0$ and $u_n=1$.
The problem will arise as in Section \ref{sect2.3} when we want to 
specialize   
$t=\infty$. One can follow closely the argument in Section \ref{sect2.3}
to examine the coefficients of the $E_\l$.
We will just state the
corresponding result in this case. Recall that 
$w_\l$ is the unique minimal length element of $W$ such that
$w_\l\cdot 0=\l$, $\{\a^{(i)}\}$ and $\{\l_{(i)}\}$ elements defined as in 
equations (\ref{1}) and (\ref{1b}).

\medskip
\noi  
Define the re-normalization of $E_\l(q,t,u)$ to be 
$$
\prod_{i=1}^l(1-{\bf q}^{-(\a_{j_i},\overline \l_{(i)}) -
\d_{j_i,0}(\a_0,\overline\l_{(i)})})E_\l(q,t,u)\ .
$$
This formula (modulo a $q$ factor) is obtained by applying the recursion 
formula (\ref{recursionkoor}) successively, starting with
$1$. The 
powers of the $t_i$ appearing in the expansion of this re-normalization
after the substitution $u=(t_0,1)$ 
are all {\sl negative} and the normalizing factor tends to $1$ when $t$
approaches infinity.
Therefore our desired specialization at $t=\infty$
is well defined. We denote by  $E_\l(q,\infty)$ this specialization.
Note that the coefficient of $e^\l$
in $E_\l(q,\infty)$ equals $1$. 
For each anti-dominant weight $\l$ we write $\R^\l(\infty)$ 
for the linear subspace 
spanned by  $\{E_\mu(q,\infty)\ |\ \mu\in\W\hspace{-0.1cm}\l \}$.
The polynomial $P_\l(q,\infty)$ is defined to be the  
unique 
$\W$- invariant element in $\R^\l(\infty)$ for which the coefficient of 
$e^\l$ equals
$1$.
\bigskip
\section{The representation-theoretical interpretation}
\bigskip
In this section we 
make no more reference to reducibility of the root system in 
question, but depending on the case we use the notation
$E_\l(q,\infty)$ to refer to the specialized versions of 
the nonsymmetric Macdonald polynomials or nonsymmetric Koornwinder polynomials.
\subsection{Proof of the Theorem \ref{teorema1}} The strategy is  
to study the degeneration of the recursion formulas 
(\ref{recursion}) and (\ref{recursionkoor}) for the polynomials 
$E_\l(q,\infty)$ and then to relate them with the Demazure character formula
(Theorem \ref{T1.9}). The crucial remark is that
we are only interested in the action of the operators
$G_{i,\l}(q,t)$ on the re-normalization of the $E_\l(q,t)$ when 
$(\l+\L_0,\a_i)\geq 0$. We see, after an examination of the operator 
$G_{i,\l}(q,t)$ in this situation, that the powers of $t_i$ appearing in
the description of its action are negative or zero. Because the same is true
for the re-normalization of  $E_\l(q,t)$ we can first make the specialization
at $t=\infty$. Moreover, the operators $G_{i,\l}(q,\infty)$ do not depend on
$\l$ anymore. In fact,  $G_{i,\l}(q,\infty)$ coincide
with the Demazure operators $\D_i$.
We are ready to state the following
\begin{Thm}\label{special-recursion}
Let $\l$ be a weight such that $(\l+\L_0,\a_i)\geq 0$. Then
\begin{equation}\label{specialized-recursion}
\D_iE_\l(q,\infty)=q^{-\left(\L_0,s_i\<\l\>\right)} E_{s_i\cdot\l}(q,\infty)\ .
\end{equation}
\end{Thm}
\begin{proof}
The statement is obvious for $(\l+\L_0,\a_i)=0$. Now, 
we know from Lemma \ref{lemma1.3} that if $(\l+\L_0,\a_i)>0$ we have
$$
l(w_{s_i\cdot\l})=l(s_iw_\l)=l(w_\l)+1.
$$
Therefore, if $w_\l=s_{j_p}\cdots s_{j_1}$ is a reduced decomposition 
$w_{s_i\cdot\l}=s_is_{j_p}\cdots s_{j_1}$ is also reduced. Henceforth, 
using the definition of $E_\l(q,\infty)$ and $E_{s_i\cdot\l}(q,\infty)$ 
and the 
recursion formulas (\ref{recursion}), (\ref{recursionkoor}) our conclusion 
follows.
\end{proof}
An immediate consequence of the Theorem \ref{special-recursion} is that
$$
\D_{w_\l}e^{\tilde\l}=q^{-\left(\L_0,w_\l\<\tilde\l\>\right)} 
E_{\l}(q,\infty)\ .
$$
The Theorem \ref{teorema1} follows by comparing this formula with
the Theorem \ref{T1.9}. A simple consequence of Theorem \ref{teorema1} is that
if we expand $E_\l(q,\infty)$ in terms of monomials the coefficients 
that appear are polynomials in $q^{-1}$ 
with {\sl positive integer} coefficients.
\medskip
\subsection{Proof of the Theorem \ref{teorema2}}
Let us begin with a characterization of $P_\l(q,\infty)$. If $\l$ is 
anti-dominant, $(\l,\a_i)\leq 0$ and the Theorem \ref{special-recursion} 
together with $\D_i^2=\D_i$ shows that
$$
\D_iE_\l(q,\infty)=E_\l(q,\infty)\ .
$$
This immediately implies that $E_\l(q,\infty)$ is $\W$-invariant.
\begin{Thm}
If $\l$ is an 
anti-dominant weight then
$$
P_\l(q,\infty)=E_\l(q,\infty)\ .
$$
\end{Thm}
Now, because $P_\l(q,\infty)$ is essentially the character of the Demazure
module $E_{w_\l}(\tilde\l)$ the $\W$-invariance of $P_\l(q,\infty)$ 
translates into saying that $E_{w_\l}(\tilde\l)$ decomposes into a direct sum
of simple ${\ring {\mf g}}$-modules. Let us write
$$
E_{w_\l}(\tilde\l)=
\bigoplus_{j\geq 0}E_{w_\l}(\tilde\l)_{j}
$$
where $E_{w_\l}(\tilde\l)_{j}$ is the direct sum of weight
spaces whose weights are of the form 
$\mu+j\d+(\L_0, w_\l\<\tilde \l\>)\d$ 
with integer $j$ and $\mu\in P$. Since 
$\d$ is $\W$-invariant each of the $E_{w_\l}(\tilde\l)_{j}$
decomposes as a direct sum of simple 
${\ring {\mf g}}$-modules. 
If $\chi_\mu$ is the character of $V_\mu$ the irreducible 
${\ring {\mf g}}$-module with highest weight $\mu$
$$
\chi(E_{w_\l}(\tilde\l)_{j})=
q^{-\left(\L_0, w_\l\<\tilde \l\>\right)-j}\sum_{\mu}c_{\l\mu}^j\chi_\mu.
$$
Here $c_{\l\mu}^j$ is the multiplicity of $V_\mu$ in 
$E_{w_\l}(\tilde\l)_{j}$. Summing up we find the polynomials
in $q^{-1}$ with {\sl positive integer} coefficients such that
$$
P_\l(q,\infty)=
\sum_{\mu\leq\l}d_{\l\mu}(q)\chi_\mu.
$$
The restriction on the sum comes from the triangular properties of 
$P_\l$. Let us remark that 
the positive integer numbers $d_{\l\mu}(1)$
are the multiplicities of the irreducible 
${\ring {\mf g}}$-modules in the
Demazure module $E_{w_\l}(\tilde\l)$. Also, $d_{\l\l}(q)=1$.
\medskip
\subsection{Proof of the Theorem \ref{teorema3}}
On one hand, because the coefficients of the expansion of 
$E_\l(q,\infty)$ in terms of monomials 
are polynomials in $q^{-1}$ 
with positive integer coefficients their limit at $q\to\infty$ exists. We will
denote by 
$$
E_\l(\infty,\infty)=\lim_{q\to \infty}E_\l(q,\infty).
$$
On the other hand, using Theorem \ref{teorema1} we can see that
$$
E_\l(\infty,\infty)=\chi(E_{w_\l}(\tilde\l)_{0})
$$
where $E_{w_\l}(\tilde\l)_{0}$ is the direct sum of weight
spaces whose weights are of the form 
$\mu+(\L_0, w_\l\<\tilde \l\>)\d$ 
with $\mu\in P$. It can be easily seen that
$E_{w_\l}(\tilde\l)_{0}$ is a $\bo$-module. 
Our conclusion follows if we prove that $E_{w_\l}(\tilde\l)_{0}$ is 
isomorphic to $\E_{\w_\l w_\circ}(\l_+)$ as $\bo$-modules.
As explained in the proof
of the Theorem \ref{teorema2} the vector space $E_{w_{\l_-}}(\tilde\l)_{0}$
is  also a $\ring {\mf g}$-module.
\begin{Thm}
The $\ring {\mf g}$-module $E_{w_{\l_-}}(\tilde\l)_{0}$ is the 
irreducible representation of $\ring {\mf g}$ with highest weight $\l_+$.
Furthermore, the $\bo$-modules $E_{w_\l}(\tilde\l)_{0}$ and 
$\E_{\w_\l w_\circ}(\l_+)$ are isomorphic.
\end{Thm}
\begin{proof}
By the Theorem \ref{teorema2}
we know that  the irreducible representation of 
$\ring {\mf g}$ with highest weight $\l_+$ occurs in the decomposition of 
$E_{w_{\l_-}}(\tilde\l)_{0}$ with multiplicity one.
Let us denote by $\ring V$ the copy of the irreducible representation of 
$\ring {\mf g}$ with highest weight $\l_+$ embedded in 
$E_{w_{\l_-}}(\tilde\l)_{0}$
and by $V$ the irreducible representation of ${\mf g}$
with highest weight $\L=\tilde \l+\L_0$. 
 It is easy to see that 
$E_{w_{\l_-}}(\tilde\l)_{0}$ is the $\bo$-module generated by
the weight space $V_{w(\L)}$ , where 
$w=\w_{\l}^{-1}w_\l$. From the fact that 
the space $V_{w(\L)}$ is one dimensional
and from 
$$
w(\L)=\l_- +(\L_0, w_\l\<\tilde \l\>)
$$
we deduce that $V_{w(\L)}$ is the lowest weight space of $\ring V$, and 
therefore 
$$
\ring V=E_{w_{\l_-}}(\tilde\l)_{0},
$$
both being equal with the $\bo$-module generated by
the weight space $V_{w(\L)}$. By the same argument the $\bo$-module 
$E_{w_\l}(\tilde\l)_{0}$ is generated by the one dimensional weight space
$$
V_{w_\l(\L)}={\ring V}_{\w_\l w_{\circ}(\l_+)}
$$
which also generates $\E_{\w_\l w_\circ}(\l_+)$ as a $\bo$-module. 
Our conclusion follows.
\end{proof}
The proof of Theorem \ref{teorema3} is now complete.


\end{document}